\documentclass[11pt]{article}
\usepackage{eurosym}
\usepackage{amssymb}
\usepackage{amsfonts}
\usepackage{amsmath}
\usepackage{graphicx}
\usepackage{hyperref}
\usepackage{subcaption}
\usepackage{tikz}
\usepackage{float}
\usepackage{cite}
\usepackage{float}

\newtheorem{theorem}{Theorem}

\newtheorem{corollary}[theorem]{Corollary}

\newtheorem{lemma}[theorem]{Lemma}

\newtheorem{proposition}[theorem]{Proposition}

\begin{document}

\author{George Avalos \\
Department of Mathematics\\
University of Nebraska-Lincoln \and Pelin G. Geredeli \\
Department of Mathematics\\
Iowa State University \and Boris Muha \\
Department of Mathematics\\
University of Zagreb}
\title{Rational Decay of A Multilayered Structure-Fluid PDE System }
\maketitle

\begin{abstract}
In this work, we consider a certain multilayered (thick layer) wave--(thin
layer) wave--heat (fluid) interactive PDE system. Such coupled PDE systems
have been used in the literature to describe the blood transport process in
mammalian vascular systems. In particular, the deformations of the boundary
interface (thin layer) are described via the two dimensional elastic
equation. The present work constitutes an investigation of the extent of the
stabilizing effects of the underlying fluid dissipation -- across the
boundary interface -- upon both the thick and thin structural components.
(All three PDE components evolve on their respective geometries.) In this
regard, our main result is the derivation of uniform decay rates for
classical solutions of this multilayered PDE model. To obtain these
estimates, necessary a priori inequalities for certain static multilayered
PDE models are generated here to ultimately allow an application of a
wellknown resolvent criterion for rational decay.

\vskip.3cm \noindent \textbf{Key terms:} Fluid-Structure Interaction,
Multilayered System, Semigroup, Rational Decay
\end{abstract}


\section{ Introduction}

\noindent The fluid structure interaction (FSI) phenomena constitutes a
broad area of research with applications in variety of real world problems \cite{bod,formaggia2010cardiovascular,kaltenbacher2018mathematical}. In
particular, mammalian blood vascular walls, being composed of viscoelastic
materials, undergo large deformations due to hemodynamic forces generated
during the blood transport process. As such, there is a coupling of
respective blood flow and wall deformation dynamics. This physiological
interaction between arterial walls and blood flow plays a crucial role in
the physiology and pathophysiology of the human cardiovascular system, and
can be mathematically realized by \textit{multilayered} FSI PDE. In such FSI
modeling, the blood flow is governed by the fluid flow PDE component
(incompressible Stokes or Navier Stokes); the displacements along the
elastic vascular wall are described by the structural PDE component (e.g.,
systems of elasticity). In this regard, the multilayered FSI modeling with a
view to understanding the incidence of aneurysm caused by arterial wall
deformations during the blood transportation process has recently been a
topic of great interest, see e.g.
\cite{bukavc2017multi,scotti2005fluid} and reference within. \newline

\noindent In this paper, we consider a simplified multilayered
structure-fluid interaction (FSI) system where the coupling of the 3D fluid
(blood flow) and 3D elastic (structural vascular wall) PDE components is
realized via an additional 2D elastic system on the boundary interface. 
\newline

\noindent \textbf{The PDE Model}\newline

\noindent Let the fluid geometry $\Omega _{f}\subseteq 
\mathbb{R}
^{3}$ be Lipschitz, and the structure domain $\Omega _{s}\subseteq 
\mathbb{R}
^{3}$ be a convex \textit{polyhedron} which is strictly contained in $\Omega
_{f}$ (See Figure 1). Moreover, fluid boundary %
is decomposed via $\partial \Omega
_{f}=\Gamma _{f}\cup \Gamma _{s},$ where $\Gamma _{s}=\partial \Omega _{s}$,
and so $\Gamma _{f}\cap \Gamma _{s}=\emptyset $. Thus, $\Gamma _{s}$ is the
boundary interface between fluid geometry $\Omega _{f}$ and structure
geometry $\Omega _{s}.$ The boundary interface is further decomposed via $%
\Gamma _{s}=\bigcup\limits_{j=1}^{K}\overline{\Gamma _{j}},$ where each $%
\Gamma _{j} $ is an open polygonal domain, with $\Gamma _{i}\cap \Gamma
_{j}=\emptyset $ for $i\neq j$. In addition, for $1\leq j\leq K,$ $n_{j}$ denotes the unit normal vector which is exterior to $\partial \Gamma
_{j} $. Also, as pictured in Figure 1, $\nu (x)$ denotes the unit
outward normal with respect to $\Omega _{f}$ (and so $\nu (x)$ is inward
with respect to $\Omega _{s}).$ 
\vspace{-2mm}

\begin{center}
\includegraphics[scale=0.30]{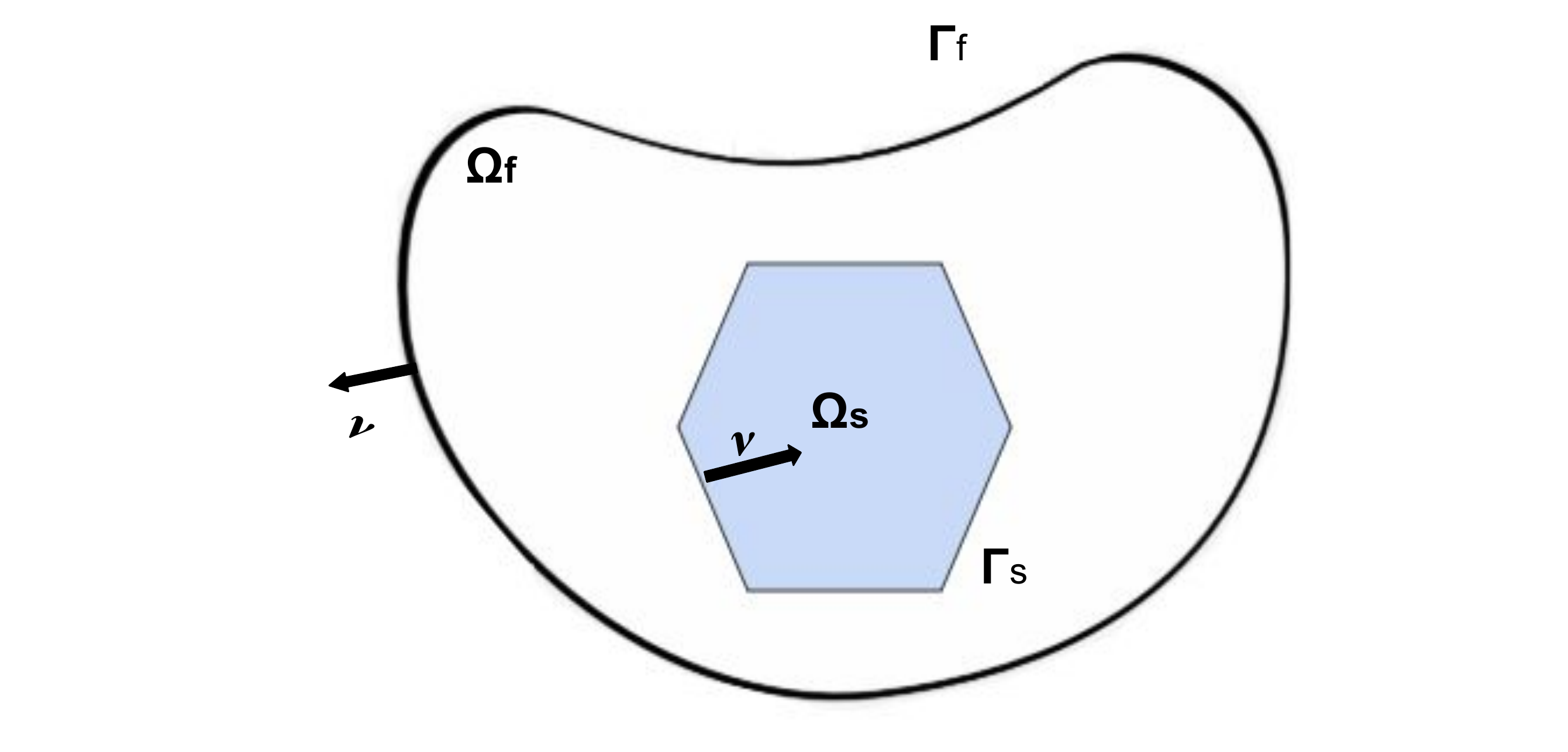}

\textbf{Figure 1: Multilayered structure-fluid interaction domain}
\end{center}

\noindent For said \{$\Omega _{f},\Gamma _{s},\Omega _{s}$\}, the
multilayered structure-fluid FSI system in solution variables $u(t,x)$ (%
corresponding to the fluid velocity), $%
h_{j}(t,x)\text{ \ }(1\leq j\leq K)$ (thin
layers displacements), and $w(t,x)$ (thick
layer displacement) is as follows:

\begin{eqnarray}
&&\left\{ 
\begin{array}{l}
u_{t}-\Delta u=0\text{ \ in }(0,T)\times \Omega _{f} \\ 
\left. u\right\vert _{\Gamma _{f}}=0\text{ \ on \ }(0,T)\times \Gamma _{f}%
\end{array}%
\right.  \label{1} \\
&&  \notag \\
&&\left\{ 
\begin{array}{l}
\text{For }1\leq j\leq K\text{,} \\ 
\frac{\partial ^{2}}{\partial t^{2}}h_{j}-\Delta h_{j}+h_{j}=\frac{\partial w%
}{\partial \nu }|_{\Gamma _{j}}-\frac{\partial u}{\partial \nu }|_{\Gamma
_{j}}\text{\ \ on }\left( 0,T\right) \times \Gamma _{j} \\ 
\\ 
\left\{ 
\begin{array}{l}
\left. h_{j}\right\vert _{\partial \Gamma _{j}\cap \partial \Gamma
_{l}}=\left. h_{l}\right\vert _{\partial \Gamma _{j}\cap \partial \Gamma
_{l}} \\ 
\left. \frac{\partial h_{j}}{\partial n_{j}}\right\vert _{\partial \Gamma
_{j}\cap \partial \Gamma _{l}}=-\left. \frac{\partial h_{l}}{\partial n_{l}}%
\right\vert _{\partial \Gamma _{j}\cap \partial \Gamma _{l}}%
\end{array}%
\right. \text{on }\left( 0,T\right) \times (\partial \Gamma _{j}\cap
\partial \Gamma _{l})\text{, \ }\forall ~~1\leq l\leq K;~~\partial \Gamma
_{j}\cap \partial \Gamma _{l}\neq \emptyset%
\end{array}%
\right.  \label{2} \\
&&  \notag \\
&&\left\{ 
\begin{array}{l}
w_{tt}-\Delta w=0\text{ \ in }(0,T)\times \Omega _{s} \\ 
\left. w_{t}\right\vert _{\Gamma _{j}}=\frac{\partial }{\partial t}%
h_{j}=\left. u\right\vert _{\Gamma _{j}}\text{ \ on }\left( 0,T\right)
\times \Gamma _{j}\text{, \ for }j=1,...,K%
\end{array}%
\right.  \label{3} \\
&&  \notag \\
&&%
\begin{array}{c}
\lbrack u(0),h_{1}(0),\frac{\partial h_{1}(0)}{\partial t},...,h_{K}(0),%
\frac{\partial h_{1}(0)}{\partial t}%
,w(0),w_{t}(0)]=[u_{0},h_{01},h_{02},...,h_{0K},h_{0K},w_{0},w_{1}]\in 
\mathbf{H}%
\end{array}
\label{4}
\end{eqnarray}%
where the finite energy space $\mathbf{H}$ is given by 
\begin{equation}
\begin{array}{l}
\mathbf{H}\equiv \left\{ \lbrack
u_{0},h_{01},h_{02},...,h_{0K},h_{0K},w_{0},w_{1}]\in \mathbf{L}^{2}(\Omega
_{f})\times \prod\limits_{j=1}^{K}\left[ H^{1}(\Gamma _{j})\times
L^{2}(\Gamma _{j})\right] \times \mathbf{H}^{1}(\Omega _{s})\times \mathbf{L}%
^{2}(\Omega _{s})\right. : \\ 
\text{(i) }\left. w_{0}\right\vert _{\Gamma _{j}}=h_{0j}; \\ 
\text{(ii) }\left. h_{j}\right\vert _{\partial \Gamma _{j}\cap \partial
\Gamma _{l}}=\left. h_{l}\right\vert _{\partial \Gamma _{j}\cap \partial
\Gamma _{l}}\text{ on }\left( 0,T\right) \times (\partial \Gamma _{j}\cap
\partial \Gamma _{l})~~\forall ~~1\leq l\leq K\text{ such that }\partial
\Gamma _{j}\cap \partial \Gamma _{l}\neq \emptyset \Big\}.%
\end{array}
\label{H}
\end{equation}%
The $\mathbf{H}$-inner product is%
\begin{eqnarray}
\left( \Phi _{0},\tilde{\Phi}_{0}\right) _{\mathbf{H}} &=&\left( u_{0,}%
\tilde{u}_{0}\right) _{\Omega _{f}}+\sum_{j=1}^{K}\left[ \left( \nabla
h_{0j},\nabla \tilde{h}_{0j}\right) _{\Gamma _{j}}+\left( h_{0j},\tilde{h}%
_{0j}\right) _{\Gamma _{j}}\right]  \notag \\
&&+\sum_{j=1}^{K}\left( h_{1j},\tilde{h}_{1j}\right) _{\Gamma _{j}}+\left(
\nabla w_{0},\nabla \tilde{w}_{0}\right) _{\Omega _{s}}+\left( w_{1},\tilde{w%
}_{1}\right) _{\Omega _{s}}  \label{inner}
\end{eqnarray}%
for $\Phi _{0}=\left[ u_{0},h_{01},h_{11}(t)....,h_{0K},h_{1K},w_{0},w_{1}%
\right] $ and $\tilde{\Phi}_{0}=\left[ \tilde{u}_{0},\tilde{h}_{01},\tilde{h}%
_{11}(t)....,\tilde{h}_{0K},\tilde{h}_{1K},\tilde{w}_{0},\tilde{w}_{1}\right]
\in \mathbf{H}.$

\noindent We should note that the boundary interface does not evolve with
time. However, it is well-accepted that if the boundary interface
displacements between structure and fluid are small relative to the scale of
the geometry, the resulting FSI models are physically relevant and reliable;
(see \cite{Gunz, Lions}). \newline

\noindent \textbf{Notation}\newline

\noindent Throughout, for a given domain $D$, the norm of corresponding
space $L^{2}(D)$ will be denoted as $||\cdot ||_{D}$ (or simply $||\cdot ||$
when the context is clear). Inner products in $L^{2}(\mathcal{O})$ or $%
\mathbf{L}^{2}(\mathcal{O})$ will be denoted by $(\cdot ,\cdot )_{\mathcal{O}%
}$, whereas inner products $L^{2}(\partial \mathcal{O})$ will be written as $%
\langle \cdot ,\cdot \rangle _{\partial \mathcal{O}}$. We will also denote
pertinent duality pairings as $\left\langle \cdot ,\cdot \right\rangle
_{X\times X^{\prime }}$ for a given Hilbert space $X$. The space $H^{s}(D)$
will denote the Sobolev space of order $s$, defined on a domain $D$; $%
H_{0}^{s}(D)$ will denote the closure of $C_{0}^{\infty }(D)$ in the $%
H^{s}(D)$-norm $\Vert \cdot \Vert _{H^{s}(D)}$. We make use of the standard
notation for the boundary trace of functions defined on $\mathcal{O}$, which
are sufficently smooth: i.e., for a scalar function $\phi \in H^{s}(\mathcal{%
O})$, $\frac{1}{2}<s<\frac{3}{2}$, $\gamma (\phi )=\phi \big|_{\partial 
\mathcal{O}},$ which is a well-defined and surjective mapping on this range
of $s$, owing to the Sobolev Trace Theorem on Lipschitz domains (see e.g., 
\cite{necas}, or Theorem 3.38 of \cite{McLean}). Also, $C>0$ will denote a
generic constant.

\section{ Literature}

Stability analysis of fluid structure interaction (FSI) PDE systems have
been an ongoing object of study \cite{16,A-T-2,mmas,A-P2,Duy,RZZ}. Because
of their utility in mathematically describing fluid or flow dynamics as they
interact with elastic materials, such FSI
models arise in biomedicine, biomechanics and aeroelasticity, %
see e.g. \cite{bod,dow}. The main
motivation of the current problem comes from the mathematical modeling of
vascular blood flow: the corresponding modeling PDE dynamics accounts for
the fact that blood-transporting vessels are generally composed of several
layers. Such multilayered FSI PDE models
have a crucial role in understanding the physiology of the human
cardiovascular system \cite{mc1, mc2,SunBorMulti}.

\noindent Examples of single layered FSI -- i.e., only one elastic PDE 
(describing three
dimensional bulk elasticity or some lower-dimensional model of plate/shell type)
models the structural dynamics; the displacement along the interaction
interface is not modeled via any elastic equation-- appear extensively in
the literature,see e.g.
\cite{AT,Barbu,benevsova2020variational,Chambolle,CS05,IKLT17,Gunz,SunBorFSIKoiter} and references within%
. However, many biomedical devices (such as stents) are being developed with
the view that vascular wall structures are composed of composite materials
and not of single layer; see \cite{multi-layered,CGLjMTW-stent,CGM-stent,
SunBorMulti, muha1}. In short, some degree of physical realism is lost if
the FSI PDE does not adequately describe arterial wall layers of composite
type.

\noindent Compared to the extensive work undertaken for single layered FSI,
there is a relative paucity of results for multilayered FSI systems. A
multilayered FSI (2D heat-1D wave- 2D wave) system was initially studied in 
\cite{SunBorMulti} with a focus on showing wellposedness. Therein, the
authors exploited an underlying regularity which was available by the
presence of the additional wave equation. A
simplified $1D$ model was studied in \cite{BorisSimplifiedFSI} where the
optimal regularity result was proved. In \cite{AGM}, wellposedness and
strong stability of a higher dimensional %
linear version of the system considered in 
\cite{SunBorMulti} were studied. In particular, to prove strong decay, the
authors of \cite{AGM} appealed to the wellknown spectral criteria in \cite%
{A-B}. Very recently, an alternative \textit{resolvent criterion approach}
to strong decay, with respect to a multilayered Lam$\acute{e}$-heat system
was given in \cite{Av2}. However, up to the present time, there has not
been, to the best of our knowledge, any investigations into \emph{uniform
decay properties }-- with respect to either finite energy or higher norm --
of multilayered FSI PDE models. Accordingly, the question in the present
work is, \textit{is the dissipation emanating from the thermal component of
the FSI system (\ref{1})-(\ref{4}) strong enough to elicit polynomial decay?}

\section{ Novelty and Challenges}

As noted above, despite extensive research activity on single-layer FSI
models in the last twenty years or so, a comprehensive long term analysis
theory for multilayered FSI -- in which the boundary interface coupling
between fluid and structure components is realized via an additional elastic
equation -- is largely absent. Having established in \cite{AGM} the strong
(asymptotic) stability for the multilayered FSI model (\ref{1})-(\ref{4}),
the authors in the present work address the issue of obtaining rational
decay rates for solutions of the coupled 3D heat-2D wave-3D wave dynamics
under consideration. These rational decay rates will pertain to solutions of
(\ref{1})-(\ref{4}) which correspond to smooth initial data; i.e., initial
data drawn from the domain of the associated heat-wave-wave $C_{0}$%
-semigroup generator $\mathbf{A}:D(\mathbf{A})\subset \mathbf{H}\rightarrow 
\mathbf{H}$ in (\ref{A-op}) below.

\noindent Our approach here for solving the rational decay problem will
entail an appropriate estimation of the resolvent of the corresponding
semigroup generator, with a view of invoking the wellknown resolvent
criterion in \cite{B-T} for polynomial decay. Ultimately, we will obtain an
explicit decay rate of $\mathcal{O}(t^{-\frac{2}{11}})$. This is Theorem \ref%
{stable} below. As as far as we can tell, this is the first such polynomial
stability result obtained for multilayered FSI. By way of obtaining the
rational decay result, we will operate in the frequency domain and deal with
a static FSI system -- which is essentially the \textit{Laplace transformed}
version of the original system (\ref{1})-(\ref{4})-- and the resolvent of
the generator of the dynamical system. In this regard, challenging issues associated
with the analysis are as follows:

\vspace{0.2cm}\noindent \textbf{(i) Majorizing the solution to the resolvent
equation in terms of the heat component.} Having obtained in our previous
work \cite{AGM} an understanding of the spectral properties of the
corresponding semigroup generator $\mathbf{A}:D(\mathbf{A})\subset \mathbf{H}%
\rightarrow \mathbf{H}$ (of (\ref{A-op}) below), we proceed with considering
the resolvent system in \eqref{stat} below and discern an inherent (static)
fluid dissipative relation. Analogous to what the analysis undertaken in the
time domain for control of PDE's in general -- see e.g., \cite{trigg} -- we
will strive here to exploit the static dissipation by \emph{majorizing} the
solution $\Phi $ of (\ref{resolve}) in terms of the heat component $\mathbf{u%
}$. However the \emph{key} issue here will be an appropriate estimate for
the 3D \textit{thick} wave component $w$.

\vspace{0.2cm}\noindent \textbf{(ii) Sharpening Poincar$\acute{e}$'s
Inequality for the thermal component.} In order to deal with critical
boundary trace terms and ultimately majorize them with respect to the static
heat dissipation, we will need to refine some of the estimates concerning
the boundary term $\left. \frac{\partial u}{\partial \nu }\right\vert
_{\Gamma _{s}}$. This will require us to prove a sharpening of Poincar$%
\acute{e}$'s Inequality; in particular, we will need to properly majorize
the $L^{2}$-norm of the thermal component in such a way so as to ultimately
secure the given decay rate.

\vspace{0.2cm}\noindent \textbf{(iii) Control of the critical three wave
boundary terms.} As we pointed out in (i), the main challenge here will be
the appropriate estimate for the 3D \textit{thick} wave component $w.$ The
use of a certain \textquotedblleft Dirichlet\textquotedblright\ map will
ultimately enable us to homogenize this thick wave component of (\ref{stat}%
), via a new variable $z$, which has zero Dirichlet bounday trace. To this
new variable $z$, we subsequently apply frequency domain versions of known
vector identities for the control of (uncoupled) waves. However, the
obtained preliminary estimate on the thick wave component will still contain
the problematic boundary term $\left. \frac{\partial z}{\partial \nu }%
\right\vert _{\Gamma _{s}}$, a term which should be controlled in $L^{2}$%
-sense. The desired estimate of this flux term will require the invocation
of the thin layer $h$-equations in (\ref{stat}) and some \textit{sharp}
interpolation inequalities. It is the estimation of this boundary term which
ultimately dictates the obtained rational decay rate.

\section{Preliminaries}

It is shown in \cite{AGM} that the multilayered PDE system (\ref{1})-(\ref{4}%
) can be described via the matrix operator $\mathbf{A}:D(\mathbf{A})\subset 
\mathbf{H}\rightarrow \mathbf{H}$ defined by 
\begin{equation}
\mathbf{A}=\left[ 
\begin{array}{cccccccc}
\Delta & 0 & 0 & \ldots & 0 & 0 & 0 & 0 \\ 
0 & 0 & I & \ldots & 0 & 0 & 0 & 0 \\ 
-\left. \frac{\partial }{\partial \nu }\right\vert _{\Gamma _{1}} & (\Delta
-I) & 0 & \ldots & 0 & 0 & \left. \frac{\partial }{\partial \nu }\right\vert
_{\Gamma _{1}} & 0 \\ 
\vdots & \vdots & \vdots & \ldots & \vdots & \vdots & \vdots & \vdots \\ 
0 & 0 & 0 & \ldots & 0 & I & 0 & 0 \\ 
-\left. \frac{\partial }{\partial \nu }\right\vert _{\Gamma _{K}} & 0 & 0 & 
\ldots & (\Delta -I) & 0 & \left. \frac{\partial }{\partial \nu }\right\vert
_{\Gamma _{K}} & 0 \\ 
0 & 0 & 0 & \ldots & 0 & 0 & 0 & I \\ 
0 & 0 & 0 & \ldots & 0 & 0 & \Delta & 0%
\end{array}%
\right] ,  \label{A-op}
\end{equation}
with

$D(\mathbf{A})=\left\{ \left[
u_{0},h_{01},h_{11}(t)....,h_{0K},h_{1K},w_{0},w_{1}\right] \in \mathbf{H:}%
\right. $

\begin{description}
\item[(\textbf{A.i)}] (a) $u_{0}\in H^{1}(\Omega _{f})$, ~(b) $h_{1j}\in
H^{1}(\Gamma _{j})$ for $j=1,...,K$, ~(c) $w_{1}\in H^{1}(\Omega _{f});$

\item[\textbf{(A.ii)}] (a) $\Delta u_{0}\in L^{2}(\Omega _{f})$, (b) $\Delta
w_{0}\in L^{2}(\Omega _{s})$, (c) $\Delta h_{0j}-\left. \frac{\partial u_{0}%
}{\partial \nu }\right\vert _{\Gamma _{j}}+\left. \frac{\partial w_{0}}{%
\partial \nu }\right\vert _{\Gamma _{j}}\in L^{2}(\Gamma _{j})$ \ for $%
j=1,...,K,$

~~(d) $\left. \frac{\partial h_{0j}}{\partial n_{j}}\right\vert _{\partial
\Gamma _{j}}\in H^{-\frac{1}{2}}(\partial \Gamma _{j})$ \ for $j=1,...,K;$

\item[(A.iii)] (a) $\left. u_{0}\right\vert _{\Gamma _{f}}=0$,~~ (b) $\left.
u_{0}\right\vert _{\Gamma j}=h_{1j}=\left. w_{1}\right\vert _{\Gamma j}$ \
for $j=1,...,K;$

\item[(A.iv)] For $j=1,...,K$ such that $\partial \Gamma _{j}\cap \partial
\Gamma _{l}\neq \emptyset $:

$\left. \text{(a)}\left. h_{1j}\right\vert _{\partial \Gamma _{j}\cap
\partial \Gamma _{l}}=\left. h_{1l}\right\vert _{\partial \Gamma _{j}\cap
\partial \Gamma _{l}};~~~(b)\left. \frac{\partial h_{0j}}{\partial n_{j}}%
\right\vert _{\partial \Gamma _{j}\cap \partial \Gamma _{l}}=-\left. \frac{%
\partial h_{0l}}{\partial n_{l}}\right\vert _{\partial \Gamma _{j}\cap
\partial \Gamma _{l}}\right\} .$
\end{description}

This is to say, $\Phi (t)=\left[ u(t),h_{1}(t),\frac{\partial }{\partial t}%
h_{1}(t)....,h_{K}(t),\frac{\partial }{\partial t}h_{K},w(t),w_{t}(t)\right]$
satisfies the PDE model\newline
(\ref{1})-(\ref{4}) if and only if these variables solve the following ODE
in Hilbert space $\mathbf{H}:$ 
\begin{equation}
\frac{d}{dt}\Phi (t)=\mathbf{A}\Phi (t)~~ \text{on}~ (0,T);~~~~~\Phi
(0)=\Phi _{0}=[u_{0},h_{01},h_{02},...,h_{0K},h_{0K},w_{0},w_{1}]\in \mathbf{%
H}.  \label{solute}
\end{equation}

\noindent We recall the wellposedness result given in \cite{AGM}:

\begin{theorem}
\label{wp} The linear operator $\mathbf{A}:D(\mathbf{A})\subset \mathbf{H}%
\rightarrow \mathbf{H},$ as defined in \eqref{A-op}, generates a $C_{0}$%
-semigroup $\left\{ e^{\mathbf{A}t}\right\} _{t\geq 0}$ of contractions on $%
\mathbf{H}.$ Thus, for $\Phi
_{0}=[u_{0},h_{01},h_{02},...,h_{0K},h_{0K},w_{0},w_{1}]\in \mathbf{H},$ the
solution $\Phi (t)=\left[ u(t),h_{1}(t),\frac{\partial }{\partial t}%
h_{1}(t)....,h_{K}(t),\frac{\partial }{\partial t}h_{K},w(t),w_{t}(t)\right]$
of (\ref{1})-(\ref{4}) is given (continuously) by 
\begin{equation*}
\Phi (t)=e^{\mathbf{A}t}\Phi _{0}\in C([0,T];\mathbf{H}).
\end{equation*}
\end{theorem}

\section{ Main Result}

\noindent Our present work mainly focuses on analyzing the long time
behavior of solutions to the given multilayered system (\ref{1})-(\ref{4}),
with a view of obtaining rational decay rate of these solutions. Our main
proof of stability will be based on an ultimate appeal to wellknown
resolvent criterion of A. Borichev and Y. Tomilov \cite[Theorem 2.4]{B-T}:

\begin{theorem}
\label{BT} Let $\left\{ T(t)\right\} _{t\geq 0}$ be a bounded $C_{0}$%
-semigroup on a Hilbert space $H$ with generator $A$ such that $i\mathbb{R}%
\subset \rho (A)$. Then for fixed $\alpha >0$ the following are equivalent:%
\vspace{0.2cm}

(i) $\left\Vert \mathcal{R}(is;A)\right\Vert =\mathcal{O}(\left\vert
s\right\vert ^{\alpha })$, \ $\left\vert s\right\vert \rightarrow \infty ;$

(ii) $\left\Vert T(t)A^{-1}x\right\Vert =o(t^{-\frac{1}{\alpha }})$, \ $%
t\rightarrow \infty $, $x\in \mathbf{H}.$
\end{theorem}
Given this operator theoretic result, it will suffice to establish the
\textquotedblleft frequency domain\textquotedblright\ PDE estimate in %
\eqref{decay}; the proof of this estimate will constitute the bulk of the
effort in the present paper. Now, we give our main result of polynomial
decay for solutions which correspond to smooth initial data as follows:

\begin{theorem}
\label{stable}In regard to the multilayered PDE system in (\ref{1})-(\ref{4}%
) (or equivalently (\ref{solute})): if $\Phi
_{0}=[u_{0},h_{01},h_{02},...,h_{0K},h_{0K},w_{0},w_{1}]\in D(\mathbf{A})$,
then the corresponding solution of (\ref{1})-(\ref{4}) (or equivalently (\ref%
{solute})) satisfies the estimate,%
\begin{equation}
\left\Vert \Phi (t)\right\Vert _{\mathbf{H}}\leq \frac{C}{t^{\frac{2}{11}}}%
\left\Vert \Phi _{0}\right\Vert _{D(\mathbf{A})}.  \label{decay}
\end{equation}%
That is, the solution to (\ref{1})-(\ref{4}) (or equivalently (\ref{solute}%
)), which corresponds to smooth initial data, decays at a rate of $\mathcal{O%
}(t^{-\frac{2}{11}})$.
\end{theorem}

\vspace{0.2cm}

\noindent \textbf{Proof of Theorem \ref{stable}}\newline

\noindent The proof relies on the resolvent criterion given in Theorem \ref%
{BT}, and presupposes that there is no intersection of $\sigma (\mathbf{A})$
with the imaginary axis. In fact it was shown in \cite[Proposition 7, Lemma
9, Corollary 10]{AGM} -- See Section 4 therein -- that $i\mathbb{R}\subset
\rho (\mathbf{A})$.

\noindent Subsequently, given parameter $\beta \in 
\mathbb{R}
$ and data $\Phi _{0}^{\ast }=[u^{\ast },h_{01}^{\ast },h_{02}^{\ast
},...,h_{0K}^{\ast },h_{0K}^{\ast },w_{0}^{\ast },w_{1}^{\ast }]\in \mathbf{H%
}$, we consider the resolvent equation%
\begin{equation}
\left[ i\beta I-\mathbf{A}\right] \Phi =\Phi _{0}^{\ast },  \label{resolve}
\end{equation}%
with solution $\Phi =[u,h_{01},h_{02},...,h_{0K},h_{0K},w_{0},w_{1}]\in D(%
\mathbf{A})$. From the definition of $\mathbf{A,}$ this abstract equation
can be written explicitly as

\begin{equation}
\begin{array}{l}
\left\{ 
\begin{array}{l}
i\beta u-\Delta u=u^{\ast }\text{ \ \ in }\Omega _{f} \\ 
u|_{\Gamma _{f}}=0\text{ on }\Gamma _{f}%
\end{array}%
\right. \\ 
\\ 
\left\{ 
\begin{array}{l}
i\beta h_{0j}-h_{1j}=h_{0j}^{\ast }\text{ \ \ in \ }\Gamma _{j} \\ 
-\beta ^{2}h_{0j}-\Delta h_{0j}+h_{0j}+\frac{\partial u}{\partial \nu }-%
\frac{\partial w_{0}}{\partial \nu }=h_{1j}^{\ast }+i\beta h_{0j}^{\ast }%
\text{ \ \ in \ }\Gamma _{j}%
\end{array}%
\right. \\ 
\left\{ 
\begin{array}{l}
\text{For all }1\leq l\leq K\text{\ such that }\partial \Gamma _{j}\cap
\partial \Gamma _{l}\neq \emptyset \\ 
h_{0j}|_{\partial \Gamma _{j}\cap \partial \Gamma _{l}}=h_{0l}|_{\partial
\Gamma _{j}\cap \partial \Gamma _{l}}\text{ } \\ 
\left. \frac{\partial h_{0j}}{\partial n_{j}}\right\vert _{\partial \Gamma
_{j}\cap \partial \Gamma _{l}}=-\left. \frac{\partial h_{0l}}{\partial n_{l}}%
\right\vert _{\partial \Gamma _{j}\cap \partial \Gamma _{l}},%
\end{array}%
\right. \\ 
\\ 
\left\{ 
\begin{array}{l}
w_{1}=i\beta w_{0}-w_{0}^{\ast }\text{ \ in }\Omega _{s}\text{\ } \\ 
-\beta ^{2}w_{0}-\Delta w_{0}=i\beta w_{0}^{\ast }+w_{1}^{\ast }\text{ \ \
in }\Omega _{s} \\ 
\lbrack i\beta w_{0}-w_{0}^{\ast }]_{\Gamma _{j}}=h_{1j}=u|_{\Gamma _{j}}%
\text{ \ on \ }\Gamma _{j}.%
\end{array}%
\right.%
\end{array}
\label{stat}
\end{equation}

\vspace{0.2cm}\noindent We will give our proof step-wise, estimating each
solution component separately:\vspace{0.2cm}

\noindent \textbf{Step 1: A static dissipation relation for heat component $%
u $: } \vspace{0.2cm}

\noindent We will start with an inherent (static) fluid dissipative relation
which will be the key ingredient for future steps. First, we take the $%
\mathbf{H}$-inner product of both sides of (\ref{resolve}) with respect to
pre-image $\Phi $. This gives%
\begin{equation}
\begin{array}{l}
i\beta \left\Vert \Phi \right\Vert _{\mathbf{H}}^{2}+\left\Vert \nabla
u\right\Vert _{\Omega _{f}}^{2}-2i\sum_{j=1}^{K}\left[ \text{Im}\left(
\nabla h_{1j},\nabla h_{0j}\right) _{\Gamma _{j}}+\text{Im}\left(
h_{1j},h_{0j}\right) _{\Gamma _{j}}\right]  \\ 
\text{ \ \ }-2i\text{Im}\left( \nabla w_{1},\nabla w_{0}\right) _{\Omega
_{s}}=\left( \Phi _{0}^{\ast },\Phi _{0}\right) _{\mathbf{H}},%
\end{array}
\label{d1}
\end{equation}%
and then the following dissipation relation: 
\begin{equation}
\left\Vert \nabla u\right\Vert _{\Omega _{f}}^{2}=\text{Re}\left( \Phi
_{0}^{\ast },\Phi \right) _{\mathbf{H}}.  \label{d2}
\end{equation}%
In view of relation (\ref{d2}), \emph{we should strive to} \emph{majorize } 
\emph{solution} $\Phi $ \emph{of} (\ref{resolve}) \emph{in norm} \emph{in
terms of the static heat dissipation}. With this theme in mind, from the
mechanical compatibility conditions in (\ref{H}), and the resolvent
relations and matching velocity BC's in (\ref{stat}), we have for $j=1,...,K$%
,%
\begin{equation}
\left[ i\beta h_{0j}-h_{0j}^{\ast }\right] _{\Gamma _{j}}=\left.
h_{1j}\right\vert _{\Gamma _{j}}=\left. u\right\vert _{\Gamma _{j}}.
\label{d3}
\end{equation}%
Combining this relation with the Sobolev Trace Theorem (and Poincar\'{e}'s
Inequality), we then have for $j=1,...,K$,%
\begin{equation}
\left\Vert \beta h_{0j}\right\Vert _{H^{\frac{1}{2}}(\Gamma _{j})}\leq
C\left( \left\Vert \nabla u\right\Vert _{\Omega _{f}}+\left\Vert \Phi
_{0}^{\ast }\right\Vert _{\mathbf{H}}\right) .  \label{d4}
\end{equation}%
Moreover, via an integration by parts we get 
\begin{eqnarray}
\left\Vert \frac{\partial u}{\partial \nu }\right\Vert _{H^{-\frac{1}{2}%
}\left( \partial \Omega _{f}\right) } &\leq &C\left( \left\Vert \Delta
u\right\Vert _{\Omega _{f}}+\left\Vert \nabla u\right\Vert _{\Omega
_{f}}\right)   \notag \\
&=&C\left( \left\Vert i\beta u-u^{\ast }\right\Vert _{\Omega
_{f}}+\left\Vert \nabla u\right\Vert _{\Omega _{f}}\right)   \label{d45}
\end{eqnarray}%
which, in turn, gives 
\begin{equation}
\left\Vert \frac{\partial u}{\partial \nu }\right\Vert _{H^{-\frac{1}{2}%
}\left( \partial \Omega _{f}\right) }\leq C\left( \left\vert \beta
\right\vert \left\Vert u\right\Vert _{\Omega _{f}}+\left\Vert \nabla
u\right\Vert _{\Omega _{f}}+\left\Vert \Phi _{0}^{\ast }\right\Vert _{%
\mathbf{H}}\right) .  \label{d5}
\end{equation}%
We should note that this estimate can be refined with respect to $\beta $.
In fact, for our particular multilayered PDE model, we have the following
\textquotedblleft sharpening\textquotedblright\ of Poincar\'{e}'s
Inequality: \ 

\begin{proposition}
\label{refine}For $\left\vert \beta \right\vert >0$, the heat solution
component of (\ref{stat}) obeys the estimate%
\begin{equation}
\left\vert \beta \right\vert ^{\frac{1}{2}}\left\Vert u\right\Vert _{\Omega
_{f}}\leq C\left( \left\Vert \nabla u\right\Vert _{\Omega _{f}}+\left\Vert
\Phi _{0}^{\ast }\right\Vert _{\mathbf{H}}\right).  \label{sharp}
\end{equation}
\end{proposition}

\noindent \textbf{Proof:} The details of the proof are taken in large part
from Lemma 5.2 of \cite{AT} and \cite{ALT}. Given heat component $u$ of (\ref%
{stat}), let variable $u_{1}$ solve%
\begin{equation}
\left\{ 
\begin{array}{l}
\Delta u_{1}=i\beta u\text{ \ in }\Omega _{f} \\ 
\left. u_{1}\right\vert _{\partial \Omega _{f}}=0\text{ \ on }\partial
\Omega _{f}.%
\end{array}%
\right.   \label{u2}
\end{equation}%
Taking the $L^{2}$-inner product of both sides of (\ref{u2})$_{1}$ by $u_{1}$%
, we subsequently have, from the heat equation in (\ref{stat}):%
\begin{eqnarray*}
\left\Vert \nabla u_{1}\right\Vert _{\Omega _{f}}^{2} &=&-\left( \Delta
u_{1},u_{1}\right) _{\Omega _{f}} \\
&=&-\left( i\beta u,u_{1}\right) _{\Omega _{f}} \\
&=&-\left( \Delta u+u^{\ast },u_{1}\right) _{\Omega _{f}} \\
&=&\left( \nabla u,\nabla u_{1}\right) _{\Omega _{f}}-\left( u^{\ast
},u_{1}\right) _{\Omega _{f}};
\end{eqnarray*}%
whence we obtain, via Poincar\'{e}'s and Young's Inequalities,%
\begin{equation}
\left\Vert \nabla u_{1}\right\Vert _{\Omega _{f}}\leq C\left( \left\Vert
\nabla u\right\Vert _{\Omega _{f}}+\left\Vert \Phi _{0}^{\ast }\right\Vert _{%
\mathbf{H}}\right) .  \label{u3}
\end{equation}%
Since the $u_{1}$-equation in (\ref{u2}) gives%
\begin{equation*}
\left\Vert u\right\Vert _{H^{-1}(\Omega _{f})}=\frac{1}{\left\vert \beta
\right\vert }\left\Vert \Delta u_{1}\right\Vert _{_{H^{-1}(\Omega _{f})}},
\end{equation*}%
using again the Poincar\'{e}'s Inequality yields%
\begin{equation}
\left\Vert u\right\Vert _{H^{-1}(\Omega _{f})}\leq \frac{C}{\left\vert \beta
\right\vert }\left\Vert \nabla u_{1}\right\Vert _{_{\Omega _{f}}}\leq \frac{C%
}{\left\vert \beta \right\vert }\left( \left\Vert \nabla u\right\Vert
_{\Omega _{f}}+\left\Vert \Phi _{0}^{\ast }\right\Vert _{\mathbf{H}}\right) ,
\label{u4}
\end{equation}%
after using (\ref{u3}). Moreover, it is known-- see e.g., Theorem B.8,
Theorem 3.30 and Theorem 3.33 of \cite{McLean}-- that 
\begin{equation*}
L^{2}(\Omega _{f})=\left[ H^{-1}(\Omega _{f}),H^{1}(\Omega _{f})\right] _{%
\frac{1}{2}}
\end{equation*}%
(which is the Lipschitz domain version of the interpolation result Lemma
12.1 of \cite{L-M}). Combining this with the estimate (\ref{u4}) gives%
\begin{equation*}
\left\Vert u\right\Vert _{L^{2}(\Omega _{f})}\leq C\left\Vert u\right\Vert
_{H^{-1}(\Omega _{f})}^{\frac{1}{2}}\left\Vert \nabla u\right\Vert _{\Omega
_{f}}^{\frac{1}{2}}\leq \frac{C}{\left\vert \beta \right\vert ^{\frac{1}{2}}}%
\left( \left\Vert \nabla u\right\Vert _{\Omega _{f}}+\left\Vert \Phi
_{0}^{\ast }\right\Vert _{\mathbf{H}}\right) 
\end{equation*}%
which yields (\ref{sharp}) and finishes the proof of Proposition \ref{refine}%
. \ \ \ $\square $ \vspace{0.2cm}

\noindent Now, applying the estimate (\ref{sharp}) to the right hand side of
(\ref{d5}), we then have the following normal derivative trace estimate for
the heat component of the static problem \eqref{stat}:

\begin{corollary}
\label{normal}For $\left\vert \beta \right\vert >1$, the heat solution
component of (\ref{stat}) obeys the estimate%
\begin{equation}
\left\Vert \frac{\partial u}{\partial \nu }\right\Vert _{H^{-\frac{1}{2}%
}\left( \partial \Omega _{f}\right) }\leq C\left\vert \beta \right\vert ^{%
\frac{1}{2}}\left( \left\Vert \nabla u\right\Vert _{\Omega _{f}}+\left\Vert
\Phi _{0}^{\ast }\right\Vert _{\mathbf{H}}\right) .  \label{u5}
\end{equation}
\end{corollary}

\vspace{0.2cm}

\noindent \textbf{Step 2: The thick wave displacement $w_0$} \vspace{0.3cm}

\noindent In what follows, we require \textquotedblleft
Dirichlet\textquotedblright\ map $D$, defined by having for given boundary
function $g\in H^{\frac{1}{2}}(\Gamma _{s})$,%
\begin{equation}
\Delta Dg=0\text{ \ in }\Omega _{s};\ \ \ \ \ \left. Dg\right\vert _{\Gamma
_{s}}=g\text{ \ on }\Gamma _{s}.  \label{D}
\end{equation}%
By the Lax-Milgram Theorem and an argument similar to that which resulted in
(\ref{d5}), we have%
\begin{equation}
D\in \mathcal{L}\left( H^{\frac{1}{2}}(\Gamma _{s}),H^{1}(\Omega
_{s})\right) \text{, \ }\frac{\partial D}{\partial \nu }\in \mathcal{L}%
\left( H^{\frac{1}{2}}(\Gamma _{s}),H^{-\frac{1}{2}}(\Gamma _{s})\right) .
\label{D2}
\end{equation}%
(The latter is the \textquotedblleft Dirichlet to Neumann\textquotedblright\
map.) Therewith, and with respect to the thick wave displacement $w_{0}$ in (%
\ref{stat}), we set 
\begin{equation}
z=w_{0}+\frac{i}{\beta }D\left( \left. u\right\vert _{\Gamma _{s}}+\left.
w_{0}^{\ast }\right\vert _{\Gamma _{s}}\right) .  \label{z}
\end{equation}%
Then, via (\ref{D}), the variable $z$ satisfies the following boundary value
problem:%
\begin{eqnarray}
-\beta ^{2}z-\Delta z &=&-i\beta D\left[ \left. u\right\vert _{\Gamma
_{s}}+\left. w_{0}^{\ast }\right\vert _{\Gamma _{s}}\right] +w_{1}^{\ast
}+i\beta w_{0}^{\ast }\text{ \ in }\Omega _{s};  \notag \\
z &=&0\text{ \ on }\Gamma _{s}.  \label{z2}
\end{eqnarray}%
Since $\Omega _{s}$ is convex then $z\in H^{2}(\Omega _{s})$ (see e.g.,
Theorem 3.2.1.2, p. 147 of \cite{grisvard}). Subsequently, we can appeal to
the known Sobolev boundary regularity results for polyhedral domains; \cite[%
p. 43, Theorem 6.9]{dauge}. In short, we have the estimate, for $\left\vert
\beta \right\vert \geq 1$,%
\begin{eqnarray}
\left\Vert z\right\Vert _{H^{2}(\Omega
_{s})}+\sum\limits_{j=1}^{K}\left\Vert \frac{\partial z}{\partial \nu }%
\right\Vert _{H^{\frac{1}{2}}\left( \Gamma _{j}\right) } &\leq &C_{0}\left(
\beta ^{2}\left\Vert z\right\Vert _{\Omega _{s}}+\left\vert \beta
\right\vert \left\Vert D\left[ \left. u\right\vert _{\Gamma _{s}}+\left.
w_{0}^{\ast }\right\vert _{\Gamma _{s}}\right] \right\Vert _{\Omega
_{s}}+\left\Vert w_{1}^{\ast }+i\beta w_{0}^{\ast }\right\Vert _{\Omega
_{s}}\right)   \notag \\
&\leq &C_{0}\left\vert \beta \right\vert \left( \left\Vert \beta
z\right\Vert _{\Omega _{s}}+\left\Vert \nabla u\right\Vert _{\Omega
_{f}}+\left\Vert \Phi _{0}^{\ast }\right\Vert _{\mathbf{H}}\right) ,
\label{dau}
\end{eqnarray}%
after using (\ref{z2}), (\ref{D2}) and the Sobolev Imbedding Theorem.\newline

\noindent With respect to the $z$-wave equation in (\ref{z2}), we appeal to
the \textquotedblleft frequency domain\textquotedblright\ version of the
well-known wave identity which is synonymous with boundary control of wave
equations; see Proposition 7 (ii) of \cite{mmas}, and also \cite{chen},\cite%
{trigg}. We adopt those wave identities to our solution component $z$ in the
following Proposition:

\begin{proposition}
\label{wave} Let $\mathbf{m}(x)=[m_{1}(x),m_{2}(x),m_{3}(x)]$ be an
arbitrary real-valued $[C^{2}(\bar{\Omega}_{s})]^{3}$-vector field, with
associated Jacobian matrix $M(x)$. Then the wave component of the solution
to the resolvent equation (\ref{resolve}) obeys the following relation:

\begin{eqnarray}
\text{(i) \ }\int_{\Omega _{s}}\left\vert M^{\frac{1}{2}}(x)\nabla
z\right\vert ^{2}d\Omega _{s} &=&-\text{Re}\int_{\Gamma _{s}}\frac{\partial z%
}{\partial \nu }\left( \mathbf{m}\cdot \nabla \overline{z}\right) d\Gamma
_{s}+\frac{1}{2}\int_{\Gamma _{s}}\left\vert \frac{\partial z}{\partial \nu }%
\right\vert ^{2}\mathbf{m}\cdot \mathbf{\nu }d\Gamma _{s}  \notag \\
&&+\frac{1}{2}\int_{\Omega _{s}}\left\{ \left\vert \nabla z\right\vert
^{2}-\beta ^{2}\left\vert z\right\vert ^{2}\right\} \text{div}(\mathbf{m}%
)d\Omega _{s}h  \notag \\
&&-\text{Re}\int_{\Omega _{s}}\left( i\beta D\left[ \left. u\right\vert
_{\Gamma _{s}}+\left. w_{0}^{\ast }\right\vert _{\Gamma _{s}}\right]
-w_{1}^{\ast }-i\beta w_{0}^{\ast }\right) \left( \mathbf{m}\cdot \nabla 
\overline{z}\right) d\Omega _{s}  \notag \\
&&  \notag \\
&=&-\frac{1}{2}\int_{\Gamma _{s}}\left\vert \frac{\partial z}{\partial \nu }%
\right\vert ^{2}\mathbf{m}\cdot \mathbf{\nu }d\Gamma _{s}+\frac{1}{2}%
\int_{\Omega _{s}}\left\{ \left\vert \nabla z\right\vert ^{2}-\beta
^{2}\left\vert z\right\vert ^{2}\right\} \text{div}(\mathbf{m})d\Omega _{s} 
\notag \\
&&+i\beta \text{Re}\int_{\Omega _{s}}\nabla \left( D\left[ \left.
u\right\vert _{\Gamma _{s}}+\left. w_{0}^{\ast }\right\vert _{\Gamma _{s}}%
\right] -w_{0}^{\ast }\right) \cdot \mathbf{m}\bar{z}d\Omega _{s}  \notag \\
&&+i\beta \text{Re}\int_{\Omega _{s}}\left( D\left[ \left. u\right\vert
_{\Gamma _{s}}+\left. w_{0}^{\ast }\right\vert _{\Gamma _{s}}\right]
-w_{0}^{\ast }\right) \bar{z}\text{div}(\mathbf{m})d\Omega _{s}  \notag \\
&&+\text{Re}\int_{\Omega _{s}}w_{1}^{\ast }\left( \mathbf{m}\cdot \nabla 
\overline{z}\right) d\Omega _{s}.  \label{main}
\end{eqnarray}%
(ii) If $\mathbf{\tilde{m}}(x)$ is an arbitrary real-valued $[C^{2}(\bar{%
\Omega}_{s})]^{3}$-vector field, then the wave component of the solution to
the resolvent equation (\ref{resolve}) satisfies the relation,%
\begin{eqnarray}
\int_{\Omega _{s}}\left\{ \left\vert \nabla z\right\vert ^{2}-\beta
^{2}\left\vert z\right\vert ^{2}\right\} \text{div}(\mathbf{\tilde{m}}%
)d\Omega _{s} &=&-\text{Re}\int_{\Omega _{s}}\left( i\beta D\left[ \left.
u\right\vert _{\Gamma _{s}}+\left. w_{0}^{\ast }\right\vert _{\Gamma _{s}}%
\right] -w_{1}^{\ast }-i\beta w_{0}^{\ast }\right) \overline{z}\text{div}%
\left( \mathbf{\tilde{m}}\right) d\Omega _{s}  \notag \\
&&-\text{Re}\int_{\Omega _{s}}\left[ \nabla z\cdot \nabla \text{div}(\mathbf{%
\tilde{m}})\right] \bar{z}d\Omega _{s}.  \label{main_2}
\end{eqnarray}%
(Note that the expressions (\ref{main})-(\ref{main_2}) each reflect the
fact that $z=0$ on $\Gamma _{s}$ and/or the unit normal vector $\nu (x)$ is
pointing inward with respect to solid geometry $\Omega _{s}$.)
\end{proposition}

\noindent Now, let the vector fields $\mathbf{m}$ and $\mathbf{\tilde{m}}$
in (\ref{main}) and (\ref{main_2}), respectively, be taken as 
\begin{equation}
\mathbf{m}(x)=\mathbf{\tilde{m}}(x)=x.  \label{rad}
\end{equation}%
Then, via (\ref{D2}), Young's Inequality and the Sobolev Trace Theorem, we
have for $\left\vert \beta \right\vert \geq 1$,%
\begin{equation}
\left\vert \int_{\Omega _{s}}\left\{ \left\vert \nabla z\right\vert
^{2}-\beta ^{2}\left\vert z\right\vert ^{2}\right\} d\Omega _{s}\right\vert
\leq \epsilon \beta ^{2}\left\Vert z\right\Vert _{\Omega
_{s}}^{2}+C_{\epsilon }\left( \left\Vert \nabla u\right\Vert _{\Omega
_{f}}^{2}+\left\Vert \Phi _{0}^{\ast }\right\Vert _{\mathbf{H}}^{2}\right) .
\label{main_3}
\end{equation}%
In turn, with vector field $\mathbf{m}$ as specified in (\ref{rad}),
applying (\ref{main_3}) to the right hand side of (\ref{main}), using (\ref%
{D2}), the Sobolev Trace Theorem and again Young's Inequality, we have the
following Proposition:

\begin{proposition}
\label{dis}For $\left\vert \beta \right\vert \geq 1$ and $\epsilon >0$, the
variable $z$ of (\ref{z}) and (\ref{z2}) satisfies%
\begin{equation}
\int_{\Omega _{s}}\left\vert \nabla z\right\vert ^{2}d\Omega _{s}\leq
C^{\ast }\int_{\Gamma _{s}}\left\vert \frac{\partial z}{\partial \nu }%
\right\vert ^{2}d\Gamma _{s}+\epsilon \left( \left\Vert \nabla z\right\Vert
_{\Omega _{s}}^{2}+\beta ^{2}\left\Vert z\right\Vert _{\Omega
_{s}}^{2}\right) +C_{\epsilon }\left( \left\Vert \nabla u\right\Vert
_{\Omega _{f}}^{2}+\left\Vert \Phi _{0}^{\ast }\right\Vert _{\mathbf{H}%
}^{2}\right)  \label{main_4}
\end{equation}%
(with respect to (\ref{main_3}) there has also been a rescaling of $\epsilon
>0$).
\end{proposition}

At this point, we also note in the $z$-wave relation (\ref{main}) that $%
\mathbf{m}(x)$ could be specified to be the smooth vector field of Lemma
1.5.1.9, pg. 40 of \cite{grisvard}: That is, for some $\delta >0$, $\mathbf{m%
}(x)\in \lbrack C^{\infty }(\overline{\Omega }_{s})]$ satisfies 
\begin{equation}
-\mathbf{m}(x)\cdot \nu \geq \delta \text{ \ a.e. \ on }\Gamma _{s}.
\label{m}
\end{equation}%
This gives in (\ref{main}):%
\begin{eqnarray*}
\frac{\delta }{2}\int_{\Gamma _{s}}\left\vert \frac{\partial z}{\partial \nu 
}\right\vert ^{2}d\Gamma _{s} &\leq &\left\vert \int_{\Omega _{s}}\left\vert
M^{\frac{1}{2}}(x)\nabla z\right\vert ^{2}d\Omega _{s}-\frac{1}{2}%
\int_{\Omega _{s}}\left\{ \left\vert \nabla z\right\vert ^{2}-\beta
^{2}\left\vert z\right\vert ^{2}\right\} \text{div}(\mathbf{m})d\Omega
_{s}\right. \\
&&-i\beta \text{Re}\int_{\Omega _{s}}\nabla \left( D\left[ \left.
u\right\vert _{\Gamma _{s}}+\left. w_{0}^{\ast }\right\vert _{\Gamma _{s}}%
\right] -w_{0}^{\ast }\right) \cdot \mathbf{m}\bar{z}d\Omega _{s} \\
&&\left. -i\beta \text{Re}\int_{\Omega _{s}}\left( D\left[ \left.
u\right\vert _{\Gamma _{s}}+\left. w_{0}^{\ast }\right\vert _{\Gamma _{s}}%
\right] -w_{0}^{\ast }\right) \bar{z}\text{div}(\mathbf{m})d\Omega _{s}-%
\text{Re}\int_{\Omega _{s}}w_{1}^{\ast }\left( \mathbf{m}\cdot \nabla 
\overline{z}\right) d\Omega _{s}\right\vert .
\end{eqnarray*}
Estimating right hand side of this relation by means of (\ref{D2}), Young's
Inequality and the Sobolev Trace Theorem, we have the control of the trace
term $\frac{\partial z}{\partial \nu }|_{\Gamma _{s}}$ in the following
proposition:

\begin{proposition}
The variable $z$ of (\ref{z}) and (\ref{z2}) satisfies%
\begin{equation}
\int_{\Gamma _{s}}\left\vert \frac{\partial z}{\partial \nu }\right\vert
^{2}d\Gamma _{s}\leq C\left( \left\Vert \nabla z\right\Vert _{L^{2}(\Omega
_{s})}^{2}+\beta ^{2}\left\Vert z\right\Vert _{\Omega _{s}}^{2}+\left\Vert
\nabla u\right\Vert _{\Omega _{f}}^{2}+\left\Vert \Phi _{0}^{\ast
}\right\Vert _{\mathbf{H}}^{2}\right) .  \label{main_5}
\end{equation}
\end{proposition}

We should emphasize that given the right hand side of estimate (\ref{main_4}%
), it is apparent that a useful estimate of the thick wave energy component
of (\ref{stat}) will necessitate \textquotedblleft decent\textquotedblright\
control of $\left. \frac{\partial z}{\partial \nu }\right\vert _{\Gamma _{s}}
$ in $L^{2}$-sense: that is, the positive constant $C$ in (\ref{main_5})
need not be \textquotedblleft small\textquotedblright . However, while the
estimate (\ref{main_5}) does not constitute such needed control, it will
serve as an ingredient for attaining that end.

\vspace{0.2cm}

\noindent \textbf{Step 3: The thin wave displacement $h_0$} \vspace{0.3cm}

\noindent With respect to the $h$-wave equations in (\ref{stat}), we take
the $L^{2}$-inner product with respect to $h_{0j}$, for $j=1,...,K$. This
gives

\begin{eqnarray}
-\left( \Delta h_{0j},h_{0j}\right) _{\Gamma _{j}}+\left\Vert
h_{0j}\right\Vert _{\Gamma _{j}}^{2} &=&\left( \beta
^{2}h_{0j},h_{0j}\right) _{\Gamma _{j}}+\left\langle \frac{\partial w_{0}}{%
\partial \nu },h_{0j}\right\rangle _{\Gamma _{j}}-\left\langle \frac{%
\partial u}{\partial \nu },h_{0j}\right\rangle _{\Gamma _{j}}  \notag \\
&&+\left( h_{1j}^{\ast },h_{0j}\right) _{\Gamma _{j}}+i\beta \left(
h_{0j}^{\ast },h_{0j}\right) _{\Gamma _{j}}.  \label{h1}
\end{eqnarray}%
Subsequently, we invoke the Green's Theorem (for $j=1,...,K$) to get%
\begin{eqnarray}
&&\left\Vert \nabla h_{0j}\right\Vert _{\Gamma _{j}}^{2}+\left\Vert
h_{0j}\right\Vert _{\Gamma _{j}}^{2}-\left\langle \frac{\partial h_{0j}}{%
\partial \nu },h_{0j}\right\rangle _{\partial \Gamma _{j}}  \notag \\
&=&\beta ^{2}\left\Vert h_{0j}\right\Vert _{\Gamma _{j}}^{2}+\left\langle 
\frac{\partial w_{0}}{\partial \nu },h_{0j}\right\rangle _{\Gamma
_{j}}-\left\langle \frac{\partial u}{\partial \nu },h_{0j}\right\rangle
_{\Gamma _{j}}  \notag \\
&&+\left( h_{1j}^{\ast },h_{0j}\right) _{\Gamma _{j}}+i\beta \left(
h_{0j}^{\ast },h_{0j}\right) _{\Gamma _{j}}.  \label{h2}
\end{eqnarray}%
Using the thin layer boundary conditions in (\ref{stat}), we then have upon
summation%
\begin{eqnarray*}
&&\sum\limits_{j=1}^{K}\left[ \left\Vert \nabla h_{0j}\right\Vert _{\Gamma
_{j}}^{2}+\left\Vert h_{0j}\right\Vert _{\Gamma _{j}}^{2}\right] \\
&=&\sum\limits_{j=1}^{K}\left[ \beta ^{2}\left\Vert h_{0j}\right\Vert
_{\Gamma _{j}}^{2}+\left\langle \frac{\partial w_{0}}{\partial \nu }%
,h_{0j}\right\rangle _{\Gamma _{j}}-\left\langle \frac{\partial u}{\partial
\nu },h_{0j}\right\rangle _{\Gamma _{j}}\right. \\
&&\left. +\left( h_{1j}^{\ast },h_{0j}\right) _{\Gamma _{j}}+i\beta \left(
h_{0j}^{\ast },h_{0j}\right) _{\Gamma _{j}}\right] .
\end{eqnarray*}%
Estimating right hand side by means of (\ref{u5}) and (\ref{d4}), and via
the inequality%
\begin{equation*}
\left\langle \frac{\partial u}{\partial \nu },h_{0j}\right\rangle _{\Gamma
_{j}}\leq \left( \frac{1}{\left\vert \beta \right\vert ^{\frac{1}{2}}}%
\left\Vert \frac{\partial u}{\partial \nu }\right\Vert _{H^{-\frac{1}{2}%
}(\Gamma _{j})}\right) \left( \left\vert \beta \right\vert ^{\frac{1}{2}%
}\left\Vert h_{0j}\right\Vert _{H^{\frac{1}{2}}(\Gamma _{j})}\right) ,
\end{equation*}%
we then have for $\left\vert \beta \right\vert \geq 1$ that the thin wave
solution components of (\ref{stat}) obeys the following (intermediate)
estimate: 
\begin{equation}
\sum\limits_{j=1}^{K}\left[ \left\Vert \nabla h_{0j}\right\Vert _{\Gamma
_{j}}^{2}+\left\Vert h_{0j}\right\Vert _{\Gamma _{j}}^{2}\right] \leq
\left\vert \sum\limits_{j=1}^{K}\left\langle \frac{\partial w_{0}}{\partial
\nu },h_{0j}\right\rangle _{\Gamma _{j}}\right\vert +C\left( \left\Vert
\nabla u\right\Vert _{\Omega _{f}}^{2}+\left\Vert \Phi _{0}^{\ast
}\right\Vert _{\mathbf{H}}^{2}\right) .  \label{I1}
\end{equation}

\subsection{An Appropriate Estimate for $\left. \frac{\partial z}{\partial 
\protect\nu }\right\vert _{\Gamma _{s}}$}

Using the decomposition (\ref{z}) and the thin wave equations in (\ref{stat}%
), we have for $j=1,...,K,$%
\begin{equation}
\frac{\partial z}{\partial \nu }=-\beta ^{2}h_{0j}-\Delta h_{0j}+h_{0j}+%
\frac{\partial u}{\partial \nu }+\frac{i}{\beta }\frac{\partial D}{\partial
\nu }\left[ \left. u\right\vert _{\Gamma _{s}}+\left. w_{0}^{\ast
}\right\vert _{\Gamma _{s}}\right] -h_{1j}^{\ast }-i\beta h_{0j}^{\ast }%
\text{ \ \ in \ }\Gamma _{j}.  \label{norm}
\end{equation}%
For the first term on right hand side, we use the matching velocities BC and 
$h$-resolvent relation in (\ref{stat}) and a wellknown trace moment
inequality -- see e.g., Theorem 1.6.6, p. 37 of \cite{scott} -- so as to have%
\begin{equation}
\left\Vert \beta ^{2}h_{0j}\right\Vert _{\Gamma _{j}}=\left\Vert \beta
\left. u\right\vert _{\Gamma _{j}}+\beta \left. h_{0j}^{\ast }\right\vert
_{\Gamma _{j}}\right\Vert _{\Gamma _{j}}\leq C\left\vert \beta \right\vert
\left( \left\Vert u\right\Vert _{\Omega _{f}}^{\frac{1}{2}}\left\Vert \nabla
u\right\Vert _{\Omega _{f}}^{\frac{1}{2}}+\left\Vert \Phi _{0}^{\ast
}\right\Vert _{\mathbf{H}}\right) .  \label{bet0}
\end{equation}%
Subsequently invoking the improvement over Poincar\'{e}'s Inequality in
Proposition \ref{refine}, we have now for $\left\vert \beta \right\vert \geq
1$,%
\begin{equation*}
\left\Vert \beta ^{2}h_{0j}\right\Vert _{\Gamma _{j}}\leq C\left( \left\vert
\beta \right\vert ^{\frac{3}{4}}\left\Vert \nabla u\right\Vert _{\Omega
_{f}}+\left\vert \beta \right\vert \left\Vert \Phi _{0}^{\ast }\right\Vert _{%
\mathbf{H}}\right) .
\end{equation*}%
Applying this estimate to the right hand side of (\ref{norm}), along with (%
\ref{I1}), (\ref{u5}), (\ref{D2}), and the Sobolev Trace Theorem, we then
have for $\left\vert \beta \right\vert \geq 1$ (upon summing over $j=1,...,K$%
),%
\begin{equation}
\sum\limits_{j=1}^{K}\left\Vert \frac{\partial z}{\partial \nu }\right\Vert
_{H^{-1}(\Gamma _{j})}\leq C_{1}\left( \sqrt{\sum\limits_{j=1}^{K}\left\vert
\left\langle \frac{\partial w_{0}}{\partial \nu },h_{0j}\right\rangle
_{\Gamma _{j}}\right\vert }+C\left\vert \beta \right\vert ^{\frac{3}{4}%
}\left\Vert \nabla u\right\Vert _{\Omega _{f}}+\left\vert \beta \right\vert
\left\Vert \Phi _{0}^{\ast }\right\Vert _{\mathbf{H}}\right) .  \label{neg}
\end{equation}%
To refine the right hand side: Using again the decomposition (\ref{z}), we
have for $\left\vert \beta \right\vert \geq 1$, 
\begin{eqnarray*}
\sqrt{\sum\limits_{j=1}^{K}\left\vert \left\langle \frac{\partial w_{0}}{%
\partial \nu },h_{0j}\right\rangle _{\Gamma _{j}}\right\vert } &=&\sqrt{%
\sum\limits_{j=1}^{K}\left\vert \left\langle \frac{\partial }{\partial \nu }%
\left[ z-\frac{i}{\beta }D([u+w_{0}^{\ast }]|_{\Gamma _{s}})\right]
,h_{0j}\right\rangle \right\vert _{\Gamma _{j}}} \\
&=&\sqrt{\sum\limits_{j=1}^{K}\left\vert \left\langle \frac{1}{\beta }\frac{%
\partial z}{\partial \nu },\beta h_{0j}\right\rangle -\frac{i}{\beta }%
\left\langle \frac{\partial }{\partial \nu }D([u+w_{0}^{\ast }]|_{\Gamma
_{s}}),h_{0j}\right\rangle \right\vert _{\Gamma _{j}}}.
\end{eqnarray*}%
Applying (\ref{d4}), (\ref{D2}), Sobolev Trace Theorem and the Young's
Inequality, we then obtain for $\left\vert \beta \right\vert \geq 1$,%
\begin{equation*}
\sqrt{\sum\limits_{j=1}^{K}\left\vert \left\langle \frac{\partial w_{0}}{%
\partial \nu },h_{0j}\right\rangle _{\Gamma _{j}}\right\vert }\leq \frac{%
\delta ^{\ast }}{\left\vert \beta \right\vert C_{1}}\sum\limits_{j=1}^{K}%
\left\Vert \frac{\partial z}{\partial \nu }\right\Vert _{H^{-\frac{1}{2}%
}(\Gamma _{j})}+C_{\delta ^{\ast }}\left( \left\Vert \nabla u\right\Vert
_{\Omega _{f}}+\left\Vert \Phi _{0}^{\ast }\right\Vert _{\mathbf{H}}\right) ,
\end{equation*}%
where $C_{1}>0$ is the constant in (\ref{neg}). Applying this inequality to (%
\ref{neg}), we have for $\left\vert \beta \right\vert \geq 1$,%
\begin{equation}
\sum\limits_{j=1}^{K}\left\Vert \frac{\partial z}{\partial \nu }\right\Vert
_{H^{-1}(\Gamma _{j})}\leq \frac{\delta ^{\ast }}{\left\vert \beta
\right\vert }\sum\limits_{j=1}^{K}\left\Vert \frac{\partial z}{\partial \nu }%
\right\Vert _{H^{-\frac{1}{2}}(\Gamma _{j})}+C_{2,\delta ^{\ast }}\left[
\left\vert \beta \right\vert ^{\frac{3}{4}}\left\Vert \nabla u\right\Vert
_{\Omega _{f}}+\left\vert \beta \right\vert \left\Vert \Phi _{0}^{\ast
}\right\Vert _{\mathbf{H}}\right] .  \label{neg2}
\end{equation}%
Interpolating now between (\ref{main_5}) and (\ref{neg2}) we have for $%
j=1,...,K,$%
\begin{equation}
\begin{array}{l}
\left\Vert \frac{\partial z}{\partial \nu }\right\Vert _{H^{-\frac{1}{2}%
}(\Gamma _{j})}\leq C\left\Vert \frac{\partial z}{\partial \nu }\right\Vert
_{H^{-1}(\Gamma _{j})}^{\frac{1}{2}}\left\Vert \frac{\partial z}{\partial
\nu }\right\Vert _{\Gamma _{j}}^{\frac{1}{2}} \\ 
\text{ \ }\leq C_{3}\left\{ \left( \frac{\delta ^{\ast }}{\left\vert \beta
\right\vert }\sum\limits_{j=1}^{K}\left\Vert \frac{\partial z}{\partial \nu }%
\right\Vert _{H^{-\frac{1}{2}}(\Gamma _{j})}+C_{2,\delta ^{\ast }}\left[
\left\vert \beta \right\vert ^{\frac{3}{4}}\left\Vert \nabla u\right\Vert
_{\Omega _{f}}+\left\vert \beta \right\vert \left\Vert \Phi _{0}^{\ast
}\right\Vert _{\mathbf{H}}\right] \right) ^{\frac{1}{2}}\right. \\ 
\text{ \ \ \ \ \ }\left. \times \left( \left\Vert \nabla z\right\Vert
_{\Omega _{s}}+\left\Vert \beta z\right\Vert _{\Omega _{s}}+\left\Vert
\nabla u\right\Vert _{\Omega _{f}}+\left\Vert \Phi _{0}^{\ast }\right\Vert _{%
\mathbf{H}}\right) ^{\frac{1}{2}}\right\} .%
\end{array}
\label{Int1}
\end{equation}%
Subsequently, for $\left\vert \beta \right\vert \geq 1$, with constant $%
C_{3} $ as in (\ref{Int1}) and using the relation $\left\vert ab\right\vert
\leq \frac{\delta ^{\ast }KC_{3}}{2\left\vert \beta \right\vert }a^{2}+\frac{%
\left\vert \beta \right\vert }{2\delta ^{\ast }KC_{3}}b^{2}$, we then have%
\begin{eqnarray*}
\left\Vert \frac{\partial z}{\partial \nu }\right\Vert _{H^{-\frac{1}{2}%
}(\Gamma _{j})} &\leq &\frac{1}{2K}\sum\limits_{j=1}^{K}\left\Vert \frac{%
\partial z}{\partial \nu }\right\Vert _{H^{-\frac{1}{2}}(\Gamma _{j})}+\frac{%
\delta ^{\ast }KC_{3}^{2}}{2\left\vert \beta \right\vert }\left( \left\Vert
\nabla z\right\Vert _{L^{2}(\Omega _{s})}+\left\Vert \beta z\right\Vert
_{\Omega _{s}}\right) \\
&+&C_{\delta ^{\ast }}\left[ \left\vert \beta \right\vert ^{\frac{7}{4}%
}\left\Vert \nabla u\right\Vert _{\Omega _{f}}+\left\vert \beta \right\vert
^{2}\left\Vert \Phi _{0}^{\ast }\right\Vert _{\mathbf{H}}\right] \text{, \
for }j=1,...,K.
\end{eqnarray*}%
Summing this estimate over $j$ now gives%
\begin{eqnarray*}
\sum\limits_{j=1}^{K}\left\Vert \frac{\partial z}{\partial \nu }\right\Vert
_{H^{-\frac{1}{2}}(\Gamma _{j})} &\leq &\frac{1}{2}\sum\limits_{j=1}^{K}%
\left\Vert \frac{\partial z}{\partial \nu }\right\Vert _{H^{-\frac{1}{2}%
}(\Gamma _{j})}+\frac{\delta ^{\ast }}{\left\vert \beta \right\vert }%
C_{4}\left( \left\Vert \nabla z\right\Vert _{L^{2}(\Omega _{s})}+\left\Vert
\beta z\right\Vert _{\Omega _{s}}\right) \\
&&+C_{\delta ^{\ast }}\left( \left\vert \beta \right\vert ^{\frac{7}{4}%
}\left\Vert \nabla u\right\Vert _{\Omega _{f}}+\left\vert \beta \right\vert
^{2}\left\Vert \Phi _{0}^{\ast }\right\Vert _{\mathbf{H}}\right)
\end{eqnarray*}%
(where constant $C_{4}$ is independent of $\delta ^{\ast }>0$). We have then
for $j=1,...,K$,%
\begin{equation}
\left\Vert \frac{\partial z}{\partial \nu }\right\Vert _{H^{-\frac{1}{2}%
}(\Gamma _{j})}\leq \frac{\delta ^{\ast }}{\left\vert \beta \right\vert }%
C_{5}\left( \left\Vert \nabla z\right\Vert _{\Omega _{s}}+\left\Vert \beta
z\right\Vert _{\Omega _{s}}\right) +C_{\delta ^{\ast }}\left( \left\vert
\beta \right\vert ^{\frac{7}{4}}\left\Vert \nabla u\right\Vert _{\Omega
_{f}}+\left\vert \beta \right\vert ^{2}\left\Vert \Phi _{0}^{\ast
}\right\Vert _{\mathbf{H}}\right) .  \label{Int2}
\end{equation}%
Subsequently, we interpolate between (\ref{Int2}) and (\ref{dau}) to have, 
\begin{equation*}
\begin{array}{l}
\left\Vert \frac{\partial z}{\partial \nu }\right\Vert _{\Gamma _{j}}\leq
C_{6}\left\Vert \frac{\partial z}{\partial \nu }\right\Vert _{H^{-\frac{1}{2}%
}(\Gamma _{j})}^{\frac{1}{2}}\left\Vert \frac{\partial z}{\partial \nu }%
\right\Vert _{H^{\frac{1}{2}}(\Gamma _{j})}^{\frac{1}{2}} \\ 
\\ 
\text{ \ }\leq C_{6}\left\{ \left[ \frac{\delta ^{\ast }}{\left\vert \beta
\right\vert }C_{5}\left( \left\Vert \nabla z\right\Vert _{\Omega
_{s}}+\left\Vert \beta z\right\Vert _{\Omega _{s}}\right) +C_{\delta ^{\ast
}}\left( \left\vert \beta \right\vert ^{\frac{7}{4}}\left\Vert \nabla
u\right\Vert _{\Omega _{f}}+\left\vert \beta \right\vert ^{2}\left\Vert \Phi
_{0}^{\ast }\right\Vert _{\mathbf{H}}\right) \right] \right. \\ 
\text{ \ \ \ \ \ \ \ }\left. \times C_{0}\left\vert \beta \right\vert \left(
\left\Vert \beta z\right\Vert _{\Omega _{s}}+\left\Vert \nabla u\right\Vert
_{\Omega _{f}}+\left\Vert \Phi _{0}^{\ast }\right\Vert _{\mathbf{H}}\right)
\right\} ^{\frac{1}{2}} \\ 
\\ 
\leq C_{6}\left\{ \delta ^{\ast }C_{0}C_{5}\left( \left[ \left\Vert \nabla
z\right\Vert _{\Omega _{s}}+\left\Vert \beta z\right\Vert _{\Omega
_{s}}+\left\Vert \nabla u\right\Vert _{\Omega _{f}}+\left\Vert \Phi
_{0}^{\ast }\right\Vert _{\mathbf{H}}\right] ^{2}\right) \right. \\ 
\left. +C_{\delta ^{\ast }}C_{0}\left\vert \beta \right\vert \left(
\left\vert \beta \right\vert ^{\frac{7}{4}}\left\Vert \nabla u\right\Vert
_{\Omega _{f}}+\left\vert \beta \right\vert ^{2}\left\Vert \Phi _{0}^{\ast
}\right\Vert _{\mathbf{H}}\right) \left( \left\Vert \nabla z\right\Vert
_{\Omega _{s}}+\left\Vert \beta z\right\Vert _{\Omega _{s}}+\left\Vert
\nabla u\right\Vert _{\Omega _{f}}+\left\Vert \Phi _{0}^{\ast }\right\Vert _{%
\mathbf{H}}\right) \right\} ^{\frac{1}{2}}%
\end{array}%
\end{equation*}%
Invoking once more the relation $\left\vert ab\right\vert \leq \delta ^{\ast
}a^{2}+C_{\delta ^{\ast }}b^{2}$, \ for $\delta ^{\ast }>0$, we have now%
\begin{equation*}
\left\Vert \frac{\partial z}{\partial \nu }\right\Vert _{\Gamma _{j}}\leq
\delta ^{\ast }C_{7}\left( \left\Vert \nabla z\right\Vert _{\Omega
_{s}}+\left\Vert \beta z\right\Vert _{\Omega _{s}}\right) +C_{\delta ^{\ast
}}\left( \left\vert \beta \right\vert ^{\frac{11}{4}}\left\Vert \nabla
u\right\Vert _{\Omega _{f}}+\left\vert \beta \right\vert ^{3}\left\Vert \Phi
_{0}^{\ast }\right\Vert _{\mathbf{H}}\right) .
\end{equation*}%
Summing over $j=1,...,K$, and rescaling $\delta ^{\ast }>0$, we have now the
following desired trace estimate for $\left. \frac{\partial z}{\partial \nu }%
\right\vert _{\Gamma _{s}}$:

\bigskip

\begin{lemma}
\label{nz}For $\left\vert \beta \right\vert \geq 1$ and arbitrary $\delta
^{\ast }>0$, the the variable $z$ of (\ref{z}) and (\ref{z2}) obeys the
estimate%
\begin{equation}
\left\Vert \frac{\partial z}{\partial \nu }\right\Vert _{\Gamma _{s}}\leq
\delta ^{\ast }\left( \left\Vert \nabla z\right\Vert _{\Omega
_{s}}+\left\Vert \beta z\right\Vert _{\Omega _{s}}\right) +C_{\delta ^{\ast
}}\left( \left\vert \beta \right\vert ^{\frac{11}{4}}\left\Vert \nabla
u\right\Vert _{\Omega _{f}}+\left\vert \beta \right\vert ^{3}\left\Vert \Phi
_{0}^{\ast }\right\Vert _{\mathbf{H}}\right) .  \label{crux}
\end{equation}
\end{lemma}

\vspace{0.2cm}

\noindent \textit{\textbf{Completion of the Proof of Theorem \ref{stable}}}%
\vspace{0.3cm}

\noindent Applying \eqref{crux} to right hand side of (\ref{main_4}) of
Proposition \ref{dis} (and subsequently rescaling), we have for $\left\vert
\beta \right\vert \geq 1$,%
\begin{equation}
\int_{\Omega _{s}}\left\vert \nabla z\right\vert ^{2}d\Omega _{s}\leq
\epsilon \left( \left\Vert \nabla z\right\Vert _{\Omega _{s}}^{2}+\beta
^{2}\left\Vert z\right\Vert _{\Omega _{s}}^{2}\right) +C_{\epsilon }\left(
\left\vert \beta \right\vert ^{\frac{11}{2}}\left\Vert \nabla u\right\Vert
_{\Omega _{f}}^{2}+\left\vert \beta \right\vert ^{6}\left\Vert \Phi
_{0}^{\ast }\right\Vert _{\mathbf{H}}^{2}\right) .  \label{s1}
\end{equation}%
In turn, we take vector field $\mathbf{\tilde{m}}(x)$ in (\ref{main_2}) of
Proposition \ref{wave} to satisfy $\text{div}(\mathbf{\tilde{m}})=1$.
Afterwards, we estimate this relation by means of (\ref{s1}), (\ref{D2}),
the Sobolev Trace Theorem and the Young's Inequality, $\left\vert
ab\right\vert \leq \rho a^{2}+C_{\rho }b^{2}$ ($1>\rho >0$) to have%
\begin{equation}
\int_{\Omega _{s}}\left\vert \beta z\right\vert ^{2}d\Omega _{s}\leq \frac{%
\epsilon }{\left( 1-\rho \right) }\left( \left\Vert \nabla z\right\Vert
_{L^{2}(\Omega _{s})}^{2}+\beta ^{2}\left\Vert z\right\Vert _{L^{2}(\Omega
_{s})}^{2}\right) +C_{\epsilon ,\rho }\left( \left\vert \beta \right\vert ^{%
\frac{11}{2}}\left\Vert \nabla u\right\Vert _{\Omega _{f}}^{2}+\left\vert
\beta \right\vert ^{6}\left\Vert \Phi _{0}^{\ast }\right\Vert _{\mathbf{H}%
}^{2}\right) .  \label{s2}
\end{equation}%
In sum, for $\left\vert \beta \right\vert \geq 1$, (\ref{s1}) and (\ref{s2})
give now, 
\begin{equation}
\left\Vert z\right\Vert _{H^{1}(\Omega _{s})}^{2}+\left\Vert \beta
z\right\Vert _{\Omega _{s}}^{2}+\left\Vert \frac{\partial z}{\partial \nu }%
\right\Vert _{\Gamma _{s}}^{2}\leq C\left( \left\vert \beta \right\vert ^{%
\frac{11}{2}}\left\Vert \nabla u\right\Vert _{\Omega _{f}}^{2}+\left\vert
\beta \right\vert ^{6}\left\Vert \Phi _{0}^{\ast }\right\Vert _{\mathbf{H}%
}^{2}\right) .  \label{s3}
\end{equation}%
Subsequently, via the change of variable in (\ref{z}), (\ref{s3}) and the
regularity for $D$ in (\ref{D2}), with respect to the thick layer wave
solution components in (\ref{stat}), we have that for $\left\vert \beta
\right\vert \geq 1$,%
\begin{eqnarray}
\text{(i) }\left\Vert w_{0}\right\Vert _{H^{1}(\Omega _{s})}^{2}
&=&\left\Vert z-\frac{i}{\beta }D\left( \left. u\right\vert _{\Gamma
_{s}}+\left. w_{0}^{\ast }\right\vert _{\Gamma _{s}}\right) \right\Vert
_{H^{1}(\Omega _{s})}^{2}  \notag \\
&\leq &C\left( \left\vert \beta \right\vert ^{\frac{11}{2}}\left\Vert \nabla
u\right\Vert _{\Omega _{f}}^{2}+\left\vert \beta \right\vert ^{6}\left\Vert
\Phi _{0}^{\ast }\right\Vert _{\mathbf{H}}^{2}\right) ;  \label{s4} \\
&&  \notag \\
\text{(ii) }\left\Vert w_{1}\right\Vert _{L^{2}(\Omega _{s})}^{2}
&=&\left\Vert i\beta \left[ z-\frac{i}{\beta }D\left( \left. u\right\vert
_{\Gamma _{s}}+\left. w_{0}^{\ast }\right\vert _{\Gamma _{s}}\right) \right]
-w_{0}^{\ast }\right\Vert _{L^{2}(\Omega _{s})}^{2}  \notag \\
&\leq &C\left( \left\vert \beta \right\vert ^{\frac{11}{2}}\left\Vert \nabla
u\right\Vert _{\Omega _{f}}^{2}+\left\vert \beta \right\vert ^{6}\left\Vert
\Phi _{0}^{\ast }\right\Vert _{\mathbf{H}}^{2}\right) .  \label{s5}
\end{eqnarray}%
In turn, via (\ref{I1}) and the decomposition (\ref{z}) we get%
\begin{eqnarray*}
\sum\limits_{j=1}^{K}\left[ \left\Vert \nabla h_{0j}\right\Vert _{\Gamma
_{j}}^{2}+\left\Vert h_{0j}\right\Vert _{\Gamma _{j}}^{2}\right] &\leq
&\left\vert \sum\limits_{j=1}^{K}\left\langle \frac{\partial (z-\frac{i}{%
\beta }D\left( \left. u\right\vert _{\Gamma _{s}}+\left. w_{0}^{\ast
}\right\vert _{\Gamma _{s}}\right) }{\partial \nu },h_{0j}\right\rangle
_{\Gamma _{j}}\right\vert +C\left( \left\Vert \nabla u\right\Vert _{\Omega
_{f}}^{2}+\left\Vert \Phi _{0}^{\ast }\right\Vert _{\mathbf{H}}^{2}\right) \\
&\leq &C\left( \left\Vert \frac{\partial z}{\partial \nu }\right\Vert
_{\Gamma _{s}}+\frac{1}{\left\vert \beta \right\vert }\left\Vert \frac{%
\partial \frac{i}{\beta }D\left( \left. u\right\vert _{\Gamma _{s}}+\left.
w_{0}^{\ast }\right\vert _{\Gamma _{s}}\right) }{\partial \nu }\right\Vert
_{H^{-\frac{1}{2}}(\Gamma _{s})}\right) \sum\limits_{j=1}^{K}\left\Vert
h_{0j}\right\Vert _{H^{\frac{1}{2}}(\Gamma _{j})}.
\end{eqnarray*}%
Invoking now, (\ref{s3}), (\ref{D2}), the Sobolev Trace Theorem, and (\ref%
{d4}), we have for $\left\vert \beta \right\vert \geq 1$,%
\begin{equation}
\sum\limits_{j=1}^{K}\left[ \left\Vert \nabla h_{0j}\right\Vert _{\Gamma
_{j}}^{2}+\left\Vert h_{0j}\right\Vert _{\Gamma _{j}}^{2}\right] \leq
C\left( \left\vert \beta \right\vert ^{\frac{11}{2}}\left\Vert \nabla
u\right\Vert _{\Omega _{f}}^{2}+\left\vert \beta \right\vert ^{6}\left\Vert
\Phi _{0}^{\ast }\right\Vert _{\mathbf{H}}^{2}\right) .  \label{s6}
\end{equation}%
In addition, via the resolvent relation in (\ref{stat}), we have for $%
j=1,...,K$,%
\begin{equation*}
\left\Vert h_{1j}\right\Vert _{H^{\frac{1}{2}}(\Gamma _{j})}^{2}=\left\Vert
i\beta h_{0j}-h_{0j}^{\ast }\right\Vert _{H^{\frac{1}{2}}(\Gamma _{j})}^{2},
\end{equation*}%
whence by (\ref{d4}),%
\begin{equation}
\left\Vert h_{1j}\right\Vert _{H^{\frac{1}{2}}(\Gamma _{j})}^{2}\leq C\left(
\left\Vert \nabla u\right\Vert _{\Omega _{f}}^{2}+\left\Vert \Phi _{0}^{\ast
}\right\Vert _{\mathbf{H}}^{2}\right) .  \label{s7}
\end{equation}%
Finally, collecting (\ref{s4}), (\ref{s5}), (\ref{s6}) and (\ref{s7}) we
have with the solution variable \newline
$\Phi =[u,h_{01},h_{02},...,h_{0K},h_{0K},w_{0},w_{1}]$ that%
\begin{equation*}
\left\Vert \Phi \right\Vert _{\mathbf{H}}^{2}\leq C\left( \left\vert \beta
\right\vert ^{\frac{11}{2}}\left\Vert \nabla u\right\Vert _{\Omega
_{f}}^{2}+\left\vert \beta \right\vert ^{6}\left\Vert \Phi _{0}^{\ast
}\right\Vert _{\mathbf{H}}^{2}\right) .
\end{equation*}%
Invoking now the static dissipation relation (\ref{d2}) and Young's
Inequality one last time, we obtain 
\begin{eqnarray*}
\left\Vert \Phi \right\Vert _{\mathbf{H}}^{2} &\leq &C\left( \left\vert
\beta \right\vert ^{\frac{11}{2}}\left\vert \left( \Phi _{0}^{\ast },\Phi
\right) _{\mathbf{H}}\right\vert +\left\vert \beta \right\vert
^{6}\left\Vert \Phi _{0}^{\ast }\right\Vert _{\mathbf{H}}^{2}\right) \\
&\leq &\epsilon \left\Vert \Phi \right\Vert _{\mathbf{H}}^{2}+\left\vert
\beta \right\vert ^{11}\left\Vert \Phi _{0}^{\ast }\right\Vert _{\mathbf{H}%
}^{2}.
\end{eqnarray*}%
Since data $\Phi _{0}^{\ast }\in \mathbf{H}$ in (\ref{resolve}) was
arbitrary, this gives the desired resolvent bound%
\begin{equation*}
\left\Vert \mathcal{R}(i\beta ;\mathbf{A})\right\Vert _{\mathcal{L(}\mathbf{H%
})}\leq C\left\vert \beta \right\vert ^{\frac{11}{2}}.
\end{equation*}%
An appeal to Theorem \ref{BT} now concludes the proof of Theorem \ref{stable}%
.

\section{Acknowledgment}

\noindent The authors G. Avalos and Pelin G. Geredeli would like to thank
the National Science Foundation, and acknowledge their partial funding from
NSF Grant DMS-1907823.

\noindent The author Boris Muha would like to thank the Croatian Science
Foundation (Hrvatska Zaklada za Znanost), and acknowledge their partial
funding from grant number IP-2018-01-3706.

\end{document}